\documentclass[a4paper]{amsart} % ,9pt
\usepackage{amsmath}
\usepackage{amssymb}%This gives us symbols like \square
\usepackage{verbatim}%This gives us the \begin{comment} environment.
\usepackage{amsthm}%This gives us Definition style of theorem.
\usepackage{mathrsfs}%This gives us \mathscr#
\usepackage{graphicx}
\usepackage{mathtools}
\usepackage[colorlinks,pdfdisplaydoctitle,linkcolor=blue,citecolor=blue]{hyperref}
\usepackage{caption}
\usepackage{subcaption}
\usepackage{url}
\usepackage{bbm}
\usepackage[notref,notcite,color]{showkeys}
\usepackage{epstopdf}
\usepackage{pgf,tikz}
\usepackage{extarrows}
\usetikzlibrary{arrows}
\usepackage{xcolor}
\usepackage{algorithm,algpseudocode} 
\usepackage{amssymb,amsthm}    % amssymb package contains more mathematical symbols
\usepackage{graphicx}          % graphicx package enables you to paste in graphics
\usepackage{amsmath}         % best maths package on offer
\usepackage{caption}
\usepackage{subcaption}
\usepackage{ dsfont }
\usepackage{listings}
\usepackage{bbm}
\usepackage{color}
\usepackage{ dsfont }
\usepackage{enumerate}
\usepackage{todonotes}

\usepackage{pgfplots}
\pgfplotsset{compat =1.18}
\usepackage{pgfplotstable}

\renewcommand{\Delta}{\triangle}

\definecolor{darkblue}{rgb}{0,0,0.7}

\definecolor{darkgreen}{rgb}{0.01,0.75,0.24}

\newcommand{\dee}{{\rm d}}
\DeclareMathOperator*{\argmin}{{argmin}}

\def \Ee[#1]{\mathcal{E}^{\text{{#1}}}}
\def\R{\mathbf{R}}  \def\NN{\mathbf{N}}

\def\pa[#1,#2]{\frac{\partial {#1}}{\partial {#2}} }

\def\idom[#1,#2,#3]{\int_{#1}\hspace{1pt} {#2} \hspace{1pt} \text{d}{#3}}
\def\res[#1,#2]{\left.{#1}\right|_{#2}}

\def\var[#1,#2]{\langle \delta \mathcal{E}^{\text{{#1}}}({#2}),v\rangle}
\def\vars[#1,#2,#3]{\langle \delta^2\mathcal{E}^{\text{{#1}}}({#2})v,{#3}\rangle}
\def\vard[#1,#2,#3,#4]{\langle \delta\mathcal{E}^{\text{{#1}}}({#2})-\delta\mathcal{E}^{\text{{#3}}}({#4}),v\rangle}

\def\N{\mathbb{N}}

\newcommand{\xx}{\mathcal{X}}

\newcommand{\be}{\begin{equation}}
\newcommand{\en}{\end{equation}}
\newcommand{\ben}{\begin{equation*}}
\newcommand{\enn}{\end{equation*}}
\newcommand{\bea}{\begin{aligned}}
\newcommand{\ena}{\end{aligned}}

\def\ba#1\ena{\begin{align}#1\end{align}}
\def\ban#1\enan{\begin{align*}#1\end{align*}}

\theoremstyle{plain}
\newtheorem{thm}{Theorem}[section]
\newtheorem*{thm*}{Theorem}
\newtheorem{defn}[thm]{Definition}

\newtheorem{lem}[thm]{Lemma}

\newtheorem{proposition}[thm]{Proposition}

\newtheorem{remark}[thm]{Remark}

\numberwithin{equation}{section}

\begin{document}

\title[The Stochastic Steepest Descent Method for Robust Optimization]{The Stochastic Steepest Descent Method for Robust Optimization in Banach Spaces}

\author[N. K. Chada] {Neil K. Chada}
\address{Department of Mathematics , City University of Hong Kong, 83 Tat Chee Ave, Hong Kong}
\email{neilchada123@gmail.com}

\author[P.J.~Herbert] {Philip J. Herbert}
\address{Department of Mathematics, University of Sussex, Brighton, BN1 9RF, United Kingdom}
\email{p.herbert@sussex.ac.uk}

\subjclass{49M41, 65K15, 65C05}
\keywords{steepest descent, robust optimization, Banach spaces, stochastic optimization}

\begin{abstract}
Stochastic gradient methods have been a popular and powerful choice of optimization methods, aimed at either minimizing or maximizing functions. Their advantage lies in the fact that that one approximates the gradient as opposed to using the full Jacobian matrix. One research direction, related to this, has been on the application to infinite-dimensional problems, where one may naturally have a Hilbert space framework. However, there has been limited work done on considering this in a more general setup, such as where the natural framework is that of a Banach space. This article aims to address this by the introduction of a novel stochastic method, the stochastic steepest descent method (SSD). The SSD will follow the spirit of stochastic gradient descent, which utilizes Riesz representation to identify gradients and derivatives. Our choice for using such a method is that it naturally allows one to adopt a Banach space setting, for which recent applications have exploited the benefit of this, such as in PDE-constrained shape optimization. We provide a convergence theory related to this under mild assumptions. Furthermore, we demonstrate the performance of this method on a couple of numerical applications, namely a $p$-Laplacian and an optimal control problem. Our assumptions are verified in these applications.
\end{abstract}

%\textcolor{red}{Suggested abstract:
%This article is concerned with the minimisation of the stochastic functionals on Banach spaces.
%We are interested in developing and analysing an algorithm which we will refer to as the Stochastic Steepest Descent.
%The algorithm will follow the spirit of Stochastic Gradient descent, which utilises Riesz representation to identify gradients and derivatives.
%We demonstrate the convergence of our algorithm and conclude with numerical examples, which include a $p$-Laplace problem.}

\maketitle

\section{Introduction} 
\label{sec:introduction}

Let us consider a Banach space $\mathcal{X}$, for which we will have minimal assumptions on. We define $(\Omega,\mathcal{F},\mathbb{P})$ to be a probability space and $j \colon \mathcal{X} \times \Omega \to \R$ to be a function, where we are interested in minimizing functionals $\mathcal{J}\colon \mathcal{X} \to \R$ of the form
\begin{equation}
\label{eq:func}
	\mathcal{J}(u) := \int_{\Omega} j(\xi,u) \ \dee \mathbb{P}(\xi) =: \mathbb{E}\left[ j(\cdot,u) \right].
\end{equation}
Minimising $\mathcal{J}$ is equivalent to minimising the expected value of the random map $u \mapsto j(\cdot,u)$.
This kind of problem is known to be a robust minimization problem, which allows for the 'best outcome on average', as opposed to the deterministic case of, for a given $\xi$, minimize $u\mapsto j(\xi,u)$, which provides the best case for that that specific choice of $\xi$.
While mathematically interesting, this is clearly of interest in many practical settings.
One such example is PDE-constrained optimization problems \cite{HS10,HPU08}, or specifically PDE-constrained shape optimization, \cite{DZ11,SZ92}, where the mentioned references are in the deterministic setting.

Let us state the problem of interest:
\begin{equation}
\label{eq:prob}
	\mbox{Find } u^* \in \argmin \left\{ \mathbb{E} \left[ j(\textcolor{black}{\cdot,u}) \right] : u \in \xx \right\}.
\end{equation}
Given \eqref{eq:prob}, we will not necessarily find a solution to this problem, although, our proposed algorithm will tend to a stationary point of $\mathcal{J}$.
% , which under certain (possibly unreasonable) conditions on $\mathcal{J}$, or $j$, will be a minimiser.
The typical approach to find solutions to problems of this type is to use stochastic gradient descent (SGD) methods.
This has received a lot of attention in recent years, particularly from the machine learning community as instead of requiring full gradient information to solve \eqref{eq:prob}, one uses a stochastic approximation of gradient using a random sample. 
However we are not necessarily able to proceed in this manner due to being in a Banach space.
As a result we require a different stochastic method for the problem \eqref{eq:prob}.
Our proposed method, which we aim to understand, is the stochastic steepest descent method (SSD).
To the best of our knowledge there has been no formal work which has considered the SSD in a Banach space setting.
We do, however, mention \cite{bittar2022stochastic}, who consider a stochastic auxiliary problem in Banach spaces. However unlike that work, and further work which is included below in the literature review, our work does not require assumptions such as (i) convexity, (ii) reflexitivity and (iii) separability.
As a result this formally provides our motivation for this work, aiming to introduce some important elementary results which can be verified numerically.
Before stating our highlighting result, and other contributions, we provide a brief literature review on relatable work.
\subsection{Literature review}
The SGD has been a popular choice for a range of stochastic optimization problems arising in various applied mathematical disciplines.
Since its formulation \cite{KY03,RM51}, SGD has been heavily used in the fields of machine learning, uncertainty quantification and imaging \cite{BCN18,Gower19,TJS16}. In an infinite-dimensional setting SGD operates on a Hilbert space setting, as it uses a Riesz representation.
Some important areas for which this is exploited can be found in the aforementioned references \cite{CD07,DB16,LM22} particularly related to learning in a least-squares formulation. However, the computation of the gradient does not make sense in a Banach space $\xx$ setting as the elements of the method are within the dual space $\xx^*$.
This can be an issue, when dealing with specific Banach space examples.
A natural choice of application is shape optimization, particularly PDE-constrained shape optimization \cite{DZ11}.
Typically Hilbert spaces have been used in practical settings, which remain valid in the computational setting, however recent works \cite{DJH22,PJH23,MPR23} have utilized Banach spaces, where one makes use of a non-stochastic steepest descent-type method.
In the random shape optimization setting, there has been recent work closely related to the Hilbertian setting \cite{GLW21,GSW23}, which consider inner product spaces using the metric of an infinite dimensional Riemannian manifold.
Another related direction, is PDE-constrained optimization under uncertainty. This has been explored in the context of a Hilbert setting, through the works of Nobile et al.~\cite{martin2021complexity,MN21}.
In these works various algorithms have been proposed, which also extend to optimal control problems as well. However their work is not easily extendable to the setting of Banach spaces, assuming the same algorithms are used.
Finally in the context of inverse problems a recent works.
The works \cite{GFH24,JK23} consider the application of the SGD to inverse problems in Banach spaces.
Specifically they consider a linear setting, which assumes the use of duality maps and their randomness is induced from mini-batching. There is very limite work on the SSD, in general, for the non-function space setting within traditional optimization \cite{stacey,sarowar}.
\par
\textcolor{black}{It is important to emphasize that our motivation, in terms of presenting an alternative to SG methods (such as the SSD method), is primarily in the context of PDE settings, where one may naturally find that the function space of interest is a Banach space. Our focus is not in the general optimization procedure of minimizing functions over the set of reals, as SGD methods already work extremely well.
% In other words, its a sufficiently regular Banach space, which why they work effectively well.
Hence why we review literature related to the problems arising from PDE and tomography-based problems.} 
%\todo{Hi Neil, I got rid of this Hilbertian part, I think it is hard to write in a professional way, I think it is just better to state that the other approach works well.  I am not sure how to appropriately fit in the last sentence. Neil: no problem thanks Philip it reads well}

\subsection{Main Theorem}
Our main theorem we state below, which is a convergence theorem for the SSD method in a Banach space setting. We note that the abuse the following notation of $j$ below which should be noted is different to that of \eqref{eq:func}.

\begin{thm*}Suppose $j$ is a.s.~in $C^{1,1}_{\rm loc}(\xx)$ and the sequence $\{u_n\}_{n\geq 1}$ generated by Algorithm \ref{alg:ssd} is $\mathcal{F}_n$-measurable and bounded. %\ph{not sure why it wouldn't be}.
    Further suppose that $\{t_n\}_{n \in \N} \in \ell^2 \setminus \ell^1$,
	then if 
    \begin{equation*}
    \mathcal{K}_n:= \mathbb{E}[\mathcal{J}'(u_n)[v_n] |\mathcal{F}_n ] \leq 0,
    \end{equation*}
    the sequence $\mathcal{J}(u_n)$ converges a.s.~and
    $\liminf_{n \to \infty} (-\mathcal{K}_n) = 0$.
Furthermore, if the sequence $t_n$ satisfies 
\begin{equation}
\label{eq:t}
\sum_{j=1}^{\infty}\frac{t_j}{\sum_{n=1}^j t_n}  =\infty,
\end{equation}
then, almost surely
\begin{equation}
\label{eq:rate}
\min_{n=1,\ldots,j} (-\mathcal{K}_n) = \mathcal{O}\bigg( \Big(\sum_{n=1}^j t_n\Big)^{-1} \bigg).
\end{equation}
\end{thm*}
To complement our main theorem, we provide our highlighting contributions of this work below, which include the above theorem:
\begin{itemize}
\item We motivate and discuss the use of the stochastic steepest descent (SSD) method for optimization problems. Specifically we consider this for robust optimization problems (with the addition of uncertainty) in the context of Banach spaces.
\item We provide a number of theoretical insights justifying the use of the SSD method in a Banach space framework. 
These include verifying there is a solution to our problem, differentiability of the random functional and providing a convergence theorem, as stated through the above theorem.
\item Various applications will be considered, for which we will test the SSD method on. These will include one generic function to optimize, as a toy example, while the other two will include two partial differential equations; a $p$-Laplace type problem and an optimal control problem with $L^p$ cost, for $p \geq 1$. Within each example we will prove various assumptions that are required in the infinite-dimensional case.
\item Numerical simulations are provided related to the examples above. We compare the SSD method to SGD, and demonstrate numerical gains in terms of the error and derivatives with respect the number of iterations.
\end{itemize}
%To complement our main theorem we will provide numerical experiments for two model problems.
%For each numerical experiment we will verify the assumptions for the theorem to hold. The first being a $p$-Laplace-type problem with random data, and the second being an optimal control problem with $L^p$ cost, for $p > 1$. Our assumptions are proven in the infinite-dimensional case, and therefore \textcolor{blue}{hold in the finite-dimensional approximating space in which we may conduct experiments}.

\subsection{Outline}
The outline of this paper is as follows: in Section \ref{sec:back} we provide some necessary background material.
This includes a discussion on the stochastic gradient descent method and the stochastic steepest descent method, with various assumptions.
These will lead onto Section \ref{sec:theory} where we provide out main theorem and analysis related to the SSD in Banach spaces, specifically a convergence result.
We will then discuss two particular applications of interest, which are presented in Section \ref{sec:app}.
Finally numerical experiments are provided in Section \ref{sec:num} on the two applications, where we conclude with a conclusion section.

\section{Background material}
\label{sec:back}

In this section we will provide an overview of the relevant background material. This will include a discussion on both the SSD and the SGD and their respective differences. As well as this we will cover traditional assumptions related to both, and discuss some of their properties in a function space setting. This will include both generic Hilbert and Banach spaces.

Consider $\xx$ a Banach space and $(\Omega,\mathcal{F}, \mathbb{P})$ a probability space $\Omega$ with filtration $\mathcal{F}$. We recall firstly some definitions, before discussing our algorithm of interest.
We further denote $H$ to be a Hilbert space.
\begin{defn}[Stochastic process]
Given a Banach space $\mathcal{X}$, a discrete-time stochastic process is a collection of $\mathcal{X}$-values random variables indexed by $n$, i.e. $\{\vartheta_n:\Omega \rightarrow \mathcal{X}: n \in \NN\}$.
\end{defn}

\begin{defn}[Filtration]
A filtration is a sequence $\{\mathcal{F}_n\}$ of sub $\sigma$-algebras of $\mathcal{F}$, 
such that $\mathcal{F}_1 \subset \mathcal{F}_2 \subset \ldots \subset \mathcal{F}$.
We further state that a stochastic process is said to be adapted to a filtration $\{\mathcal{F}_n\}$
if and only if, the sequence $\{\vartheta_n\}$ is $\mathcal{F}_n$-measurable for all $n \in \NN$.
\end{defn}
\begin{defn}[Gaussian measure on $\mathcal{X}$]
\label{def:gauss}
Let $\mathcal{X}$ be a Banach space. A probability measure $\mathbb{P}$ on $(\mathcal{X}, \mathcal{B}(\mathcal{X})$
is said to be Gaussian if $\mathbb{P} \circ f^{-1}$ is a Gaussian measure on $\R$ for every $f \in \mathcal{X}^*$.
The measure $\mathbb{P}$ is called centred, or symmetric, if all the measures $\mathbb{P} \circ f^{-1}$ are centred, and called nondegenerate if for any $\mathbb{P} \neq 0$, the measure $\mathbb{P} \circ f^{-1}$ is nondegenerate.
\end{defn}

Let $j\colon \Omega \times \xx \rightarrow \overline{\R}$ be a proper functional, by that we mean there exists a value such that the functional is finite. %\todo{This is complete in my view - check optimisation and fixed point theory notes? Me too, I agree it is complete}
The aim is to find $u \in \xx$ such that the following quantity
$$
\mathbb{E}\left[ j(\xi,u) \right]= \int_{\Omega} j(\xi, u)\ \dee \mathbb{P}(\xi),
$$
is minimized. 

In the setting of a Hilbert space $H$, a typical strategy is (deterministic) gradient descent based methods.
In such a case, one seeks a $v^* \in H$ which represents the negative gradient defined as
\begin{equation}
\label{eq:neg_grad}
\langle v^*, \eta \rangle_{H} = \langle - \nabla_H \mathcal{J}(u),\eta \rangle_H := - \mathcal{J}'(u)[\eta],
\end{equation}
for all $\eta \in H$,
where $$\mathcal{J}'(u)[v] = \lim_{s \to 0}\frac{\mathcal J(u+sv)-\mathcal J(u)}{s}.$$
Equivalently, one could compute
\begin{equation}
\label{eq:min1}
v^* = \argmin_{v \in H} \Big\{\frac{1}{2}\|v\|^2_H + \mathcal{J}'(u)[v]\Big\}.
\end{equation}
Alternatively one can use the (deterministic) steepest descent method which also works 
in Banach spaces, which considers the following problem
$$
v^* \in \argmin \Big\{ \mathcal{J}'(u)[v]: v \in \mathcal{X}, \|v\|_{\mathcal{X}} \leq 1 \Big\}.
$$
The approach which we introduce is the stochastic steepest descent method (SSD), which functions as, and is defined, very similarly to the stochastic gradient descent method (SGD).
Before we discuss the SSD, let us recall the SGD and various properties and common assumptions.
We describe the SSD below in Algorithm \ref{alg:ssd}, where $j_u$ denotes the gradient of $j$ w.r.t.~the second variable $u$.
Similarly we present the SGD in Algorithm \ref{alg:sgd}.
Notice the \textit{only} difference between  Algorithms \ref{alg:ssd} and \ref{alg:sgd} is line 5.
\begin{algorithm}[h!]
\begin{algorithmic}[1]
\State \textbf{Input:} \begin{itemize}
 \item sequence of step sizes $(t_n,n\in\N)$.
 \item initial $u_0 \in \mathcal{X}$ such that $j(\cdot,u_0) < \infty$ a.s.
 \end{itemize}
 \State \textbf{Output:} 
\For {$n=1,\dots, N$}
        \State generate $\xi_n \in \Omega$ independent of previous draws.
	\State find $v_n \in \argmin\Big\{ j_u(\xi_n,u_n)[v]: \ v \in \mathcal{X},  \ \|v\|_{\mathcal{X}} \leq 1\Big\}$.
	\State iterate $u_{n+1} = u_n+t_nv_n$.
\EndFor
\end{algorithmic}
 \caption{Stochastic Steepest Descent Method (SSD)}
 \label{alg:ssd}
\end{algorithm}

\begin{algorithm}[h!]
\begin{algorithmic}[1]
\State \textbf{Input:} \begin{itemize}
 \item sequence of step sizes $(t_n,n\in\N)$.
 \item initial $u_0 \in \mathcal{X}$ such that $j(\cdot,u_0) < \infty$ a.s.
 \end{itemize}
 \State \textbf{Output:} 
\For {$n=1,\dots, N$}
        \State generate $\xi_n \in \Omega$ independent of previous draws.
	\State find $v_n \textcolor{blue}{=} \argmin_{v \in H}\Big\{ j_u(\xi_n,u_n)[v] + \frac{1}{2}\|v\|^2_H : v \in H\Big\}$.
	\State iterate $u_{n+1} = u_n+t_nv_n$.
\EndFor
\end{algorithmic}
 \caption{Stochastic Gradient Descent Method (SGD)}
 \label{alg:sgd}
\end{algorithm}
We note that the SGD necessarily has a unique minimizer in line 5, whereas SSD, one might only be able to prove that there exists a minimizer, not necessarily a unique one.
For SGD, this is of course a standard variational problem in an Hilbert space, whereas for SSD, this is the minimization within the closed unit ball, a convex linear problem.
As such, the computation of step 5 within SSD may be challenging - in the examples we present, we detail how to calculate the steepest descent.

Further notice that at each iterate, we may draw multiple samples in a so-called batch.
Let us note that one might wish to choose the size of this batch dynamically.
This is not necessary, however it may be useful for certain problems where evaluations are cheap.
Of course it is clear that we require that $j$ is a.s.~differentiable in its second variable, 
otherwise the algorithm does not make sense.
%\todo{Read over this and improve readability}
In terms of the SGD related to PDE-constrained optimization, there has been recent extensive work both numerically and analytically.
For the later, this has been concerned with deriving convergence results related to the functional \eqref{eq:func} as well as within its iterates.
Specifically a relevant result is Theorem 4.7 from \cite{GW20} which states that given a functional $\mathcal{J}$ and the optimal state $u$, with various convexity assumptions which shows
$$
\mathbb{E} \big[\|u^n-u\|_{{H}}\big] \leq \mathcal{O}\Big(\frac{1}{\sqrt{n+\nu}}\Big),
$$
where $\nu$ is some particular constant related to the strong convexity. Furthermore they also show
$$
\mathbb{E}\big[\mathcal{J}(u^n) - \mathcal{J}(u)\big] \leq  \mathcal{O}\Big(\frac{1}{n + \nu}\Big).
$$
The above results are specific to use of SGD in Hilbert spaces.
The former result is the well-known Monte Carlo rate of convergence.
The notion of strong-convexity is difficult to use in a general Banach space setting, therefore attaining such equivalent results may prove difficult. Given that it is not essential to this work, this is left as potential future work.

\section{Theory}
\label{sec:theory}
In this section we provide a number of results which will be helpful in the understanding of the SSD in Banach spaces. This corresponds to our main theorem, which is a convergence result which additionally provides an asymptotic almost sure convergence rate for the minimal values of the norm of the derivative. We provide a number of definitions for which our analysis relies on. The main assumption we have on our Banach space $\xx$ is that it should be the dual of another space, i.e.~ there is some space $Y$ such that $\xx = Y^*$.
This will be used to be able to apply the theorem of Banach–Alaoglu.
A consequence of this is that we may \emph{not} generally take $\xx$ as an $L^1(U)$ function, where we assume for now that $U \subset \R^d$, for $d \geq 1$,  is a Lipschitz domain.

\begin{defn}
\label{def:1}
    Let $\mathcal{G} \colon \xx \to \overline{\R}$.
    \begin{itemize}
    \item
    We say $\mathcal{G}$ is weak-$*$ coercive if for all $\Lambda \in \R$, the sublevel set
    $$
        \{ u \in \xx : \mathcal{G}(u) \leq \Lambda \},
    $$
    is sequentially weak-$*$ relatively compact.
    \item
    We say $\mathcal{G}$ is weak-$*$-lower-semi-continuous if for all sequences $\{u_j\}_{j\in \NN} \subset \xx$ with $u_j \overset{*}{\rightharpoonup} u$ (weak-$*$ convergent), it holds
    $$
        \mathcal{G}(u) \leq \liminf_{j\to \infty} \mathcal{G}(u_j).
    $$
    \end{itemize}
\end{defn}
%\begin{assumption}
%TBD make an assumption here on the todo notes.
%\end{assumption}
%\begin{prop}
%    Suppose that $j$ is Carath\'eodory, that is measurable in the first component and continuous in the second component.
%    Then if $u \mapsto j(\cdot,u)$ is weak-$*$-lower-semi-continuous a.s., then so is $u \mapsto \mathcal{J}(u)$.
%\end{prop}
%\begin{proof}
%    Proof of weak-$*$-lower-semi-continuity follows directly from the Dominated Convergence Theorem:\todo{Need a dominating function -- perhaps requires an assumption on $\xi \mapsto j(\xi,\cdot)$ which holds uniformly on $\xx$? Growth condition/boundedness, prehaps weak star provides the bound? check over this.}
%    \begin{equation*}
%    \begin{split}
%        \mathcal{J}(u)
%        = \int_{\Omega} j(\cdot,u) \ \dee \mathbb{P}
%        &\leq \int_{\Omega} \liminf_{j \to \infty} j(\cdot, u_j) \ \dee \mathbb{P}
%        \\
%        &\leq \liminf_{j \to \infty} \int_{\Omega} j(\cdot, u_j) \ \dee \mathbb{P} \\
%        &= \liminf_{j \to \infty} \mathcal{J}(u_j).
%    \end{split}\end{equation*}
%\end{proof}

%\begin{thm}
%    Suppose $j$ is Carath\'eodory, $u \mapsto j(\cdot,u)$ is weak-$*$-lower-semi-continuous a.s., and $\mathcal{J}$ is weak-$*$ coercive.
%    Then there is $u^* \in \xx$ such that
%    $$
%        u^* \in \argmin_{u \in \xx}  \mathcal{J}(u).
%    $$
%\end{thm}\todo{we have coercivity of $\mathcal{J}$ because we do not have a result which allows us to convert assumption on $j$}
%\begin{proof}
%This result follows from the Direct Method of the calculus of variations.
%\todo{add proof, check Rindlers book}
%\end{proof}
The next result is necessary for the proposed stochastic steepest descent method.
\begin{proposition}%\todo{Philip: Not sure if this a.s.~differentiability is ok}
    Suppose $u \mapsto j(\cdot,u)$ is differentiable a.s.~and denote this $j_u(\cdot,u) \in \xx^*$.
    Let us assume that $v \mapsto j_u(\cdot,u)[v]$ is a.s.~ weak-*-lower-semi-continuous.
    Then, there exists $v \in \xx$ such that
    \begin{equation*}
        v \in \argmin\{ j_u(\cdot,u)[\tilde{v}] : \tilde{v} \in \xx,\, \|\tilde{v}\|_\xx \leq 1\},
    \end{equation*}
    a pathwise direction of steepest descent.
\end{proposition}
\begin{proof}
    By the Banach-Alaoglu Theorem, which requires that $\xx$ is the dual space to some normed space, we know that for a particular sequence, $\{v_n\}_{n\in \NN} \subset \xx$ such that $\|v_n\|_\xx \leq 1$ and $$j_u(\cdot,u)[v_n] \to \inf \{j_u(\cdot,u)[\tilde{v}] : \tilde{v} \in \xx,\, \|\tilde{v}\|_\xx \leq 1\},$$ there is a subsequence $\{n_j\}_{j \in \NN}$ and $v^* \in \xx$ such that $\|v^*\|_\xx \leq 1$ and $v_{n_j} \overset{*}{\rightharpoonup} v^*$.
    By the assumed weak-*-lower-semi-continuity, it holds that%todo{indeed, we use sequential w*lsc}
    \begin{equation*}
        j_u(\cdot,u)[v^*]\leq \inf \{j_u(\cdot,u)[\tilde{v}] : \tilde{v} \in \xx,\, \|\tilde{v}\|_\xx \leq 1\},
    \end{equation*}
    which completes the proof.
\end{proof}
Notice that, here, the condition that $v\mapsto j_u(\cdot,u)[v]$ is weak-*-lower-semi-continuous, is a technical condition which arises from the fact we have not assumed that $\xx$ is reflexive.
This assumption has implicitly appeared in shape optimization works, see \cite[equation (12)]{PagWecFar18} for example.
%\ph{I think this is the case.  I have used this in my work, but I forget the precise reasons for the necessity.  Maybe it is fine without.  DISCUSS}

\begin{remark}
Before continuing with our theoretical results, related to our algorithm the SSD method, we would like to emphasize that the assumptions needed before, from Definition \ref{def:1}, are not required from here. Rather the assumptions that we need now are provided in Section \ref{sec:app}. We will further verify each of these assumptions for each of our test problems used in the numerical simulations later in the manuscript.
\end{remark}
%Of course it is clear that we require that $j$ is \ph{I assume a.s.} differentiable in its first variable, otherwise the algorithm does not make sense.
In order to show convergence of our method, we require a few technical results.
The first of which is a classical and important result of stochastic algorithms, which we will
not prove, but refer the reader to \cite{RS71}.

\begin{lem}[Robbins-Siegmund]\label{lem:Rob-Sie}
Let $\{\mathcal{F}_n\}_{n\geq 1}$ be a filtration.
Let $\nu_n$, $a_n$, $b_n$, and $c_n$ be non-negative random variables which are adapted to $\mathcal{F}_n$ for $n\geq 1$.
If
\begin{equation}
\begin{split}
	\mathbb{E}\left[ \nu_{n+1} | \mathcal{F}_n \right]
	\leq
	\nu_n \left( 1+ a_n\right) + b_n - c_n,\quad\mbox{and}
	\\
	\sum_{n=1}^\infty a_n < \infty, \quad \sum_{n=1}^\infty b_n < \infty, \quad \mbox{ a.s.},
\end{split}
\end{equation}
then almost surely $\nu_n$ is convergent and $\sum_{n=1}^\infty c_n< \infty$.
\end{lem}

We will say a functional $\mathcal{G} \colon \xx \to \bar{\R}$ is in $C^{1,1}_{\rm loc}(\xx)$ if it is differentiable and its derivative is locally Lipschitz, i.e.~ for any bounded set $B\subset \xx$ $\exists L_B>0$ such that $\| \mathcal{G}'(u) - \mathcal{G}'(v)\|_{\xx^*} \leq L_B \|u-v\|_{\xx}$ for any $u,\,v \in B$.
\begin{lem}
\label{lem:diff}
Suppose that $\mathcal{G} \colon \xx \to \bar{\R}$ is in $C^{1,1}_{\rm loc}(\xx)$.
Then for all $u,\,v \in \xx$ there is $L>0$, which depends on $u$ and $v$ such that
$$
    \mathcal{G}(u) - \mathcal{G}(v) \leq \mathcal{G}'(v)[u-v] + L \|u-v\|_\xx^2.
$$
\end{lem}
\begin{proof}
Let $\phi\colon \R \ni t \mapsto \mathcal{G}(u+t(v-u))$, then we have that
$$\phi'(t) = \mathcal{G}'(u+t(v-u))[v-u].$$
It also holds that $$\phi(1)- \phi(0) = \phi'(0) + \int_0^1 (\phi'(t) - \phi'(0)) \dee t.$$
As such, we have
\begin{align*}
    \mathcal{G}( v) - \mathcal{G}(u)
    =
    \mathcal{G}'(u)[v-u] + \int_0^1 ( \mathcal{G}'(u+t(v-u))[v-u] - \mathcal{G}'(u)[v-u] ) \dee t
    \\
    \leq
    \mathcal{G}(u)[v-u] + \int_0^1 L \|t ( v-u) \|_\xx  \| v-u \|_\xx \dee t,
\end{align*}
which completes the result.
\end{proof}
We are now able to provide our main result.
This result states that, under appropriate conditions, the functional along the iterations generated by Algorithm \ref{alg:ssd} converges and a function of the derivative of the functional has vanishing $\liminf$.
% This vanishing $\liminf$ condition is considered a notion of global convergence, c.f.~\cite[page 88]{HPU08}.

\begin{thm}\label{thm:main}%\todo{when is $\{u_n\}$ bounded?  Any heuristics would be nice - even if we make some assumption like convexity in the vicinity}
  %  \todo{Philip: Are these all the assumptions for $j$.  Also, did we ever get around to verifying that $u_n$ is measurable?  I know the bounded-ness was non-trivial - we don't (appear to) end up using that anyway!}
	Suppose $j$ is a.s.~in $C^{1,1}_{\rm loc}(\xx)$ and the sequences $\{u_n\}_{n\geq 1}$ and $\{v_n\}_{n\geq 1}$ generated by Algorithm \ref{alg:ssd} is $\mathcal{F}_n$-measurable and bounded. %\ph{not sure why it wouldn't be}.
    Further suppose that $\{t_n\}_{n \in \N} \in \ell^2 \setminus \ell^1$,
	then if 
    \begin{equation*}
    \mathcal{K}_n:= \mathbb{E}[\mathcal{J}'(u_n)[v_n] |\mathcal{F}_n ] \leq 0,
    % \mathcal{K}_n:= \mathbb{E}[\mathcal{J}'(u_n)[v_n] | \mathcal{F}_n] \leq 0
    % =-\|\mathcal{J}'(u_n)\|_{\xx^*}
    \end{equation*}
    the sequence $\mathcal{J}(u_n)$ converges a.s.~and
    $\liminf_{n \to \infty} \mathcal{K}_n = 0$.
    % $\liminf_{n \to \infty} \| \mathcal{J}'(u_n) \|_{\xx^*} = 0$.
Furthermore, if the sequence $t_n$ satisfies 
\begin{equation}
\label{eq:t}
\sum_{j=1}^{\infty}\frac{t_j}{\sum_{n=1}^j t_n}  =\infty,
\end{equation}
then, almost surely
\begin{equation}
\label{eq:rate}
\min_{n=1,\ldots,j} (- \mathcal{K}_n) = \mathcal{O}\bigg( \Big(\sum_{n=1}^j t_n\Big)^{-1} \bigg).
\end{equation}
\end{thm}
Before we prove this statement, let us mention that it is an open question to determine general conditions on $\mathcal{J}$ to ensure that $\mathcal{K}_n$ is non-positive.
In a stochastic gradient method, this condition is trivial, as it is the case that $\mathbb{E}[\mathcal{J}'(u)[\nabla_u j(\cdot,u)]] = - \|\mathcal{J}'(u)\|^2$.

%\todo{I disagree with the $o$-comment, I prefer $\mathcal{O}$.}
\begin{proof}
    %\todo{Did part of the proof get deleted at some point?}
	Since $\{u_n\}_{n\geq 1}$ is bounded, there is $L>0$ such that $u\mapsto j_u(\cdot,u)$ is Lipschitz with constant $L>0$ on the bounded set $B := \bigcup_{n\geq 1} \{ u_n\} \subset \xx$.
    By applying Lemma \ref{lem:diff} to $u \mapsto j(\cdot, u)$ and taking expectation, it holds that
	\begin{equation*}
		\mathcal{J}(u_{n+1}) - \mathcal{J}(u_n)
		\leq
		t_n \mathcal{J}'(u_n)[v_n] + L t_n^2 \|v_n\|_\xx^2.
	\end{equation*}
	Recall that the quantities in the above are random, depending on the draw at each step of the algorithm, also, by definition, $\|v_n\|_X \leq 1$.
	By taking expectation conditional on $\mathcal{F}_n$, one has that
	\begin{equation}
    \label{eq:expect_Lipschitz}
		\mathbb{E}[\mathcal{J}(u_{n+1}) | \mathcal{F}_n]
		-
		\mathcal{J}(u_n)
		\leq
		t_n \mathbb{E}[\mathcal{J}'(u_n)[v_n] | \mathcal{F}_n] + L t_n^2 \|v_n\|_\xx^2.
	\end{equation}
	where we have used that $u_n$ is measurable with respect to $\mathcal{F}_n$. %\ph{how would we check?}.
%	Since %\ph{this is an assumption we have not yet made} 
% we are choosing the draws at each step of the algorithm independent of the previous draws, it holds that
 % \todo{Is this necessary or relevant?}
%	\begin{equation*}
%		\mathcal{K}_n = \mathbb{E}[\mathcal{J}'(u_n)[v_n] | \mathcal{F}_n].
		% =
		% -\|\mathcal{J}'(u_n)\|_{\xx^*}.
%	\end{equation*}
% \textcolor{blue}{
% We may show that
% \begin{equation*}
%     \mathbb{E}[\mathcal{J}'(u_n)[v_n] | \mathcal{F}_n]
% 	\geq
% 	-\|\mathcal{J}'(u_n)\|_{\xx^*}
% \end{equation*}
% by Jensen's inequality.
% We are interested in the so-called \emph{angle condition} $\exists \eta \in (0,1)$ such that
% \begin{equation*}
%     \mathbb{E}[\mathcal{J}'(u_n)[v_n] | \mathcal{F}_n]
% 	\leq
% 	-\eta \|\mathcal{J}'(u_n)\|_{\xx^*}.
% \end{equation*}
% It is worth noting that, if we can only show $\mathbb{E}[\mathcal{J}'(u_n)[v_n] | \mathcal{F}_n] \leq 0$, then we can proceed with the proof, showing the convergence rate for $-\mathbb{E}[\mathcal{J}'(u_n)[v_n] | \mathcal{F}_n]$, I assume.
% }
	We therefore find ourselves in the setting of Lemma \ref{lem:Rob-Sie} with  $a_n = 0$, $b_n = t_n^2 L$,
 % $c_n = t_n \|\mathcal{J}'(u_n)\|_{\xx^*}$,
 $c_n = t_n \mathcal{K}_n$,
 and $\nu_n = \mathcal{J}(u_n)$, therefore it holds that $\mathcal{J}(u_n)$ converges a.s., and $\sum_{n = 1}^\infty t_n \|\mathcal{J}'(u_n)\|_{\xx^*}$ is finite, hence $\|\mathcal{J}'(u_n)\|_{\xx^*}$ has vanishing $\liminf$.

We now turn to the proof of \eqref{eq:rate}. We begin by defining for all $n \in \N$,
$$
\eta_n := \frac{2t_n}{\sum_{j=1}^n t_j}, \quad T_1 := - \mathcal{K}_1, \quad T_{n+1} := (1-\eta_n)T_n - \eta_n \mathcal{K}_n,
$$
where we note that $\eta_1=2$ and $\eta_n \in [0,1]$ for $n>1$ since $t_n$ is decreasing we can state
$$
0 \leq \eta_n = \frac{2t_n}{\sum_{j=1}^n t_j} \leq \frac{2t_n}{nt_n} \leq 1, \quad \mathrm{for} \ n >1.
$$
Moreover,
$$
- 2t_n \mathcal{K}_n = \sum^n_{j=1}t_jT_{n+1} + t_nT_n - \sum^{n-1}_{j=1}t_jT_n.
$$%\todo{should this have a squared on LHS?}
Then if we use this in the equation \eqref{eq:expect_Lipschitz}, we have
\begin{align*}
\mathbb{E}\Big[ \mathcal{J}(u_{n+1}) + \frac{1}{2}\sum^n_{j=1}t_jT_{n+1} | \mathcal{F}_n  \Big] &= \mathbb{E}[\mathcal{J}(u_{n+1})|\mathcal{F}_n] + \frac{1}{2}\sum^n_{j=1}t_jT_{n+1} \\
&\leq \mathcal{J}(u_n) + \frac{1}{2}\sum_{j=1}^{n-1} t_j T_n - \frac{1}{2}t_n T_n + L t_n^2.
\end{align*}
Therefore, by the Robbins-Siegmund Lemma, i.e.~ Lemma \ref{lem:Rob-Sie} with $a_n = 0$, $b_n = L t_n^2$, $c_n = \frac{1}{2}t_n T_n$, and $\nu_n = \mathcal{J}(u_n) + \sum_{j=1}^{n-1} t_j T_n$, it holds that 
$$
\Big\{ \mathcal{J}(u_n) + \frac{1}{2}\sum^{n-1}_{j=1}t_jT_n \Big\} \mbox{ converges a.s.~  and } \sum^{\infty}_{n=1}t_n T_n < \infty.
$$ 
Since $\{\mathcal{J}(u_n)\}$ converges almost surely, as demonstrated above, then we know the sequence $\{ \sum^{n-1}_{j=1}t_jT_n\}$ does so similarly. Then the fact $\sum^{\infty}_{n=1}t_n T_n < \infty$ yields
$$
\lim_{n \rightarrow \infty} \frac{t_n}{\sum^{n-1}_{j=1}t_j}\sum^{n-1}_{j=1}t_jT_n = \lim_{n \rightarrow \infty}t_nT_n=0.
$$
Using \eqref{eq:t} and the above limit, it holds that $T_n = \mathcal{O}\bigg( \Big(\sum_{j=1}^{n-1} t_j\Big)^{-1} \bigg)$ almost surely.
Using the definition of $T_n$ and that $\eta_n \in [0,1]$ one can deduce that for each $n>1$ there exists a sequence $\tilde{\eta}_j \in [0,1]$ for $j = 1,\ldots,n-1$ such that
$$
\sum^{n}_{j=1}\tilde{\eta_j} = 1,
\quad
T_n = -\sum^{n-1}_{j=1} \tilde{\eta}_j \mathcal{K}_j, %\| \mathcal{J}'(u_j) \|_{\xx^*}. 
$$
Moreover,
\begin{align*}
T_2 &= -\mathcal{K}_1, \\
T_n &\geq \sum^{n-1}_{j=1} \tilde{\eta}_j\min_{k=1,\ldots,n-1} (- \mathcal{K}_k) = \min_{k=1,\ldots,n-1} (- \mathcal{K} _k)
\geq 0, \quad \ \mathrm{for } \ n>2.
\end{align*}
Therefore the above statement for $T_n$ holds for any $n>1$, which yields the result and concludes the proof. 
%\todo{Philip: I modified this proof to be in terms of $\mathcal{K}_n$, can we check that the signs are correct and such. Neil: looks correct to me Philip}
\end{proof}
\begin{remark}
The assumed regularity of $u \mapsto {j}(\cdot,u)$ is in $C^{1,1}_{\rm loc}(\mathcal{X})$ is in line with the work of Geiersbach et al.~\cite{GLW21}.
Let us note that one may relax the assumption of boundedness of $\{u_n\}_{n\geq 1}$ if one takes the stronger assumption that $u\mapsto j(\cdot,u)$ is in $C^{1,1}(\xx)$.
\end{remark}
\begin{remark}
We note as a consequence of the Robbins-Siegmund Theorem, that Theorem \ref{thm:main} implies
that one may choose a step size of the form $t_n = n^{-1/\vartheta}$, where $\vartheta \in (0.5,1]$ as this ensures $\{t_n\}_{\N} \in \ell^2 \setminus \ell^1$.
%\todo{Is this a relevant comment, or correct?  We certainly need that $\theta > 0.5$, right?  Does the theorem suggest an adaptive step size?}
\end{remark}
% \begin{proposition}
% \label{prop:1}
% Suppose $\mathcal{J}$ is coercive, then there a.s. exists a sub-sequence $\{u_{n_k}\}_{k\geq 1}$ of the sequence generated by Algorithm \ref{alg:ssd} and weak-$*$ limit $u$ such that $u_{n_k} \rightharpoonup^* u$ in $\xx$.
% Furthermore, suppose that $\mathcal{J}$ is weak-$*$-continuous, then $\mathcal{J}(u_{n_k}) \to \mathcal{J}(u)$.
% \ph{This is different to what we discussed, I think we only know that the value of $\mathcal{J}$ decreases on average, not any 'local minimising' criteria.  It could arrive at a saddle point and be stuck, of course!}
% \end{proposition}

%\nc{Want to show uniqueness in the setting where we have convexity.}
%\\
%\ph{Need to check some conditions, this is moderately hard!  We don't actually do this in the other work, it is part of the theorem that we check!}
%\ph{ Want to assume that the sets $\{ j(u,\cdot) < c\}$ are a.s. compact, or something of this variety which would provide some notion of convergence of $u$, however this doesn't help us in shape setting!}

%\ph{Now we want to show that, should the limit $u_n$ exist in some sense, it is a stationary point}

\section{Applications}
\label{sec:app}
In this section we introduce some applications of interest, for which we will numerically test our SSD on in Section \ref{sec:num}. These will \textcolor{black}{comprise} of a elliptic partial differential equation (PDE) control problem, as well as the $p$-Laplacian (PDE). For each of these applications, we will verify the necessary assumptions related to the SSD, and describe formally how to obtain the direction of steepest descent.  

\textcolor{black}{For our applications, we will consider function spaces on an open and bounded domain $D$, which should be sufficiently that we may define Sobolev spaces.
To avoid technical details, we assume it is piecewise $C^1$.}
%\todo{maybe we should use the domain $D$, I agree it is confusing for $u$ to not be an element of $U$.}
%\begin{cases}

\subsection{Application 1}
\label{subs:app1}
As a first application, let us motivate the difference between the SGD and SSD by a deterministic example in $\R^n$.
If one wishes to find a minimum of $f \colon \R^n \to \R$, one may take the gradient descent, i.e., $$x_{n+1} = x_n - t_n\nabla f(x_n),$$ with stepsize $t_n$.
Our proposal is to instead use the steepest descent, 
$$
x_{n+1} = x_n - t_n \frac{\nabla f(x_n)}{|\nabla f(x_n)|},
$$
supposing $\nabla f(x_n) \neq 0$.

As a motivating example let us provide a simple function we are interested in minimizing, which is given as
\begin{align}
\mathcal{J}(u) &= \frac{1}{2}(u-1)^2 +  \frac{1}{2}(u+1)^2 \\ \nonumber
&= u^2+1. \label{eq:fun}
\end{align}
Considering $a \sim \mathcal{U}(\{-1,1\})$ to be the uniform random variable on two points, one may consider that $j(\xi,u) = (u-a(\xi))^2$, then one finds $\mathbb{E}[j(\xi,u)] = \mathcal{J}(u)$.

\begin{remark}
It is important to note that the example discussed is defined on the set of the reals, related to minimizing a relatively simple objective. As the domain and range are the set of reals, the example given trivially fits into a "Hilbertian-framework", implying that conventional stochastic gradient methods works well here.
Our motivation for giving this example is to show that SSD works, however we are primarily interested in more complex problems where the underlying solution is defined on an infinite dimensional Banach space.
\end{remark}
%\todo{Should we comment that 'it does not work' well - in the sense of converging to a minimiser - or leave that until later}
\begin{lem}[Finding the direction of steepest descent]
    It holds that $j_u(\xi,u)[v] = 2(u-a(\xi)) v$ for any $v \in \R$
    Furthermore, it holds that $v(u)$, the direction of steepest descent is given by $v(u) = - {\rm sgn}(u-a(\xi))$, where ${\rm sgn}$ is the sign function.
\end{lem}
%\todo{proof is trivial?}
We may also, in this setting, calculate $\mathbb{E}[ \mathcal{J}'(u)[v(u)] ]$, so as to verify the statement of Theorem \ref{thm:main}.
\begin{lem}\label{lem:OneDSteepestDescent}
    It holds that
    \begin{equation}
        \mathbb{E}[ \mathcal{J}'(u)[v(u)] ]
        =
        \begin{cases}
            2 u \quad &\mbox{if } u \leq -1,\\
            0   \quad &\mbox{if } u \in (-1,1),\\
            -2u \quad &\mbox{if } u \geq 1.
        \end{cases}
    \end{equation}
    In particular, $\mathbb{E}[ \mathcal{J}'(u)[v(u)] ] \leq 0$. 
\end{lem}
%\todo{maybe these could be combined into a single lemma, or just calcualtions?}
From this, we see that the assumptions of Theorem \ref{thm:main} are satisfied, mainly that $\mathcal{K}_n \leq 0$.

\subsection{Application 2 ($p$-Laplace)}\label{sec:app:p-laplace}
Our second application we consider a problem related to the $p$-Laplacian:
\begin{equation}
\begin{cases}
\label{eq:plap}
- \Delta_p u: = \nabla \cdot (| \nabla u |^{p-2} \nabla u) &= f, \quad \mbox{ in } U, \\
 u &=g, \quad \mbox{ on } \partial U,
\end{cases}
\end{equation}
for $p \in (1,\infty)$, where $f$ and $g$ are some sufficiently regular functions.
We set $\xx = W^{1,p}_0(U)$ with norm $\| u \|_\xx := \left( \int_U |\nabla u|^p \right)^{1/p}$.
We let $g \colon \Omega \times U \to \R$ be a given random function which satisfies {$\xi \mapsto g(\xi, \cdot) \in W^{1,p}(U)$ a.s.}.
We set the function $j$ to be given by
\begin{equation*}
    j(\xi,u) = \frac{1}{p} \int_U |\nabla(u+g(\xi,\cdot))|^p  \dee x.
\end{equation*}
The derivative of which is given by
\begin{equation}\label{eq:derivative_A1}
    j_u(\xi,u)[v] = \int_U |\nabla(u+g(\xi,\cdot))|^{p-2} \nabla (u+g(\xi,\cdot)) \cdot \nabla v\ \dee x.
\end{equation}
Finding the direction $v(u)\in \xx$ is a convex optimization problem, which may be difficult to solve.
We therefore provide the following result which is used to characterize the solution further.
This result covers the case that $j_u(\xi,u) \neq 0$ - if $j_u(\xi,u) = 0$, $v(u)$ can be chosen as any element of the unit ball of $\xx$.
\begin{lem}[Finding the direction of steepest descent]
\label{lem:direct}
    Suppose $j_u(\xi,u) \neq 0$.
    One finds that $v(u) = \mu ^{-\frac{1}{(p-1)}}\tilde{v}$, where $\mu = \left(-p^{-1} j_u(\xi,u)[\tilde{v}]\right)^{\frac{p-1}{p}}$ and $\tilde{v}\in \xx$ solves
    \begin{equation*}
        p \int_{U} |\nabla \tilde{v} |^{p-2}\nabla \tilde{v} \cdot \nabla \eta\ \dee x = - j_u(\xi, u)[\eta] \quad \forall \eta \in \xx.
    \end{equation*}
\end{lem}
\begin{proof}
%\todo{probably want a reference for KKT conditions, should find the appropriate place in Hinze's book.}
Abbreviating the dependence of $v$ on $u$, our task is to find $v \in \xx$ and $\mu \in \R$, a Lagrange multiplier, such that
\begin{equation}\label{eq:KKT-plaplace}\begin{split}
    \mu \geq 0,\quad
    \|v\|_\xx^p -1 \leq 0,\quad
    (\|v\|_\xx^p -1)\mu = 0,
    \\
    \mu \left( \|v\|_\xx^p -1 \right)' + j_u(\xi,u) = 0,
\end{split}\end{equation}
where we have taken the $p$-th power of the condition that $\|v\|_\xx \leq 1$ for convenience.
Equation \eqref{eq:KKT-plaplace} are the KKT conditions specific to \eqref{eq:plap} (page 3., \cite{HPU08}).
The first line of \eqref{eq:KKT-plaplace} corresponds to the dual feasibility, primal feasibility, and the complementary slackness, respectively; the second line may be interpreted as a stationarity condition.
For convenience, let us rewrite this final line,
\begin{equation*}
    \mu p \int_{U} |\nabla v|^{p-2}\nabla v \cdot \nabla \eta\ \dee x = - j_u(\xi, u)[\eta] \quad \forall \eta \in \xx.
\end{equation*}
In order to show $\mu >0$, we suppose $\mu = 0$, then $j_u(\xi,u)[\eta] = 0$ for all $\eta \in \xx$, which contradicts $j_u(\xi,u) \neq 0$.
Since $\mu >0$, it holds that $\|v\|_\xx = 1$ by complementary slackness.
It further holds that
\begin{equation*}
    p \int_{U} |\nabla (\mu ^{\frac{1}{(p-1)}}v)|^{p-2}\nabla (\mu ^{\frac{1}{(p-1)}}v) \cdot \nabla \eta\ \dee x = - j_u(\xi, u)[\eta] \quad \forall \eta \in \xx,
\end{equation*}
which we may solve for the unknown $\tilde{v} := \mu ^{\frac{1}{(p-1)}}v$.
Since $\|v\|_\xx = 1$, we deduce $v = \tilde{v} / \|\tilde{v}\|_\xx$, where one may compute $$p \mu^{\frac{p}{p-1}} = p \|\tilde{v}\|_\xx^p = - j_u(\xi,u)[\tilde{v}],$$
which completes the proof.
\end{proof}

\subsubsection{Verification of conditions}
%\todo{If 3.1-3.4 is removed, we may not need all of these results [(iii) and (v)] -- we might wish to put these bullet points into Application 1}
For the results we have shown, we have made assumptions on $j$, $\mathcal{J}$, and $\xx$.
These are namely that
(i) $\xx$ is the dual of a normed space,
(ii) $j$ is appropriately smooth,
%(iii) $j$ is weak-$*$-lower-semi-continuous in the second variable
(iii) $v \mapsto j_u(\cdot,u)[v]$ is weak-$*$-lower-semi-continuous.
% , and that
%(v) $\mathcal{J}$ is weak-$*$-coercive.
For simplicity, let us assume that $p\geq 2$.
%\todo{thought a paragraph would be better than an indented lemma}
\\

\paragraph{\bf Verifying condition (i)}\

Since $\xx = W^{1,p}_0(U)$ is reflexive then it is certainly the dual of a space.

\paragraph{\bf Verifying condition (ii)}\

One may see that $j$ is twice differentiable in $u$, hence its derivative is locally Lipschitz.
\begin{comment}
\todo{change the lemma to state that it is twice differentiable}
\begin{lem}
Given we have a function $j \colon \mathcal{X} \times \Omega \to \R$ which relates to the functional $\mathcal{J}:\mathcal{X} \rightarrow \R$, through \eqref{eq:func}, $j$ is sufficiently smooth, i.e.~its
first derivative $j_u$ is locally Lipschitz.
\end{lem}
\begin{proof}
It is clear that $j$ is differentiable in the second variable, the form of the derivative is given in \eqref{eq:derivative_A1}.
Its second derivative, denoted by $j_{uu}$, exists and is given by
\begin{align*}
    j_{uu}(\cdot,u)[v,w]
    &= \nabla_u \Big( \int_U |\nabla(u+g(\xi,\cdot))|^{p-2} \nabla (u+g(\xi,\cdot)) \cdot \nabla v \dee x \Big) \\
     & = \int_U ( |\nabla(u + g)|^{p-2} \nabla v \cdot \nabla w \\& \ \ \ + (p-2) |\nabla (u+g)|^{p-4} (\nabla(u+g)\cdot \nabla w) (\nabla(u+g)\cdot \nabla v) ) \dee x.
\end{align*}
Then for $p>2$, using H\"older inequalities with $$\frac{1}{p/(p-2)} + \frac{1}{p} + \frac{1}{p} =1,$$
we have
\begin{equation}\begin{split}
    j_{uu}(\cdot,u)[v,w]
    \leq&\
    (p-1)\int_U |\nabla(u+g)|^{p-2} |\nabla w| |\nabla v| \dee x \\
    \leq&\
    (p-1) \|\nabla(u+g)\|_{L^p}^{p-2} \|\nabla w\|_{L^p} \|\nabla v\|_{L^p},
\end{split}\end{equation}
which is locally bounded, hence $j_u$ is locally Lipschitz.
\end{proof}
%\todo{I forgot how I did the proof before}
%\paragraph{\bf Verifying condition (iii)}\

%Since we consider the $p$ power of the norm, it is weak-$*$-lower-semi continuous.
\end{comment}

\paragraph{\bf Verifying condition (iii)}\

By reflexivity, this holds.

%\paragraph{\bf Verifying condition (v)}\

%One may show the stronger coercivity assumption by seeing that
%\begin{equation*}
%    j(\cdot,u) \geq \frac{1}{p^2} \int_U |\nabla u|^p\dee x - p^{p-2}\int_U |\nabla g|^p \dee x.
%\end{equation*}
%After integrating in probability,
%\[
%    \mathcal{J}(u) \geq \frac{1}{p^2} \int_U |\nabla u|^p\dee x - p^{p-2} \mathbb{E}\left(\int_U |\nabla g|^p \dee x\right),
%\]
%which {assuming appropriate $p$-integrability of $g$}\todo{does this need a separate assumption?} yields the desired coercivity.

\subsection{Application 3 (optimal control problem)}\label{sec:app:optControl}
Our second application is based on a PDE constrained optimization problem, where there control-to-state map involves solving a semi-linear elliptic equation.
Given $p\in (1,\infty)$, let $\xx = L^p(U)$ for $U \subset \R^d$ for $d=2,3$ and equip it with the standard norm $\|u\|_\xx := \left( \int_U |u|^p \right)^{1/p}$.
Fix $y_d \in L^2(U)$ and $\beta>0$.
The function $j$ is given by
\begin{equation}
\label{eq:A1}
j(\xi,u) = \frac{1}{2}\int_U |y(\xi,\cdot) - y_d|^2 \dee x + \frac{\beta}{p}\int_U |u|^p \dee x,
\end{equation}
where for a random diffusivity coefficient $\kappa(\xi)$ with $0 < \kappa_0<\kappa(\xi)<\infty$ a.s., \\
% \todo{I think we might want to assume that $\kappa(\xi)\geq\kappa_0>0$ to ensure the properties we wish}
$y(\xi,\cdot)\in H^1(U)$ satisfies the following PDE
\begin{equation}
\label{eq:Poisson}
\begin{cases}
    -\nabla \cdot \left( \kappa (\xi) \nabla y \right) +y + y^5&= F(\xi) + u\quad \mbox{ in } U, \\
    \nu \cdot \nabla y &= 0\quad \ \ \ \ \  \ \ \ \  \  \mbox{ on } \partial U,
\end{cases}
\end{equation}
where again we assume $U \subset \R^d$ is a Lipschitz domain, $\partial U$ denotes the boundary of $U$ and $\nu$ is the outwards unit normal.
For $F(\xi) \in H^{-1}(U)$ and any $u \in L^p(U)$, it holds the system in \eqref{eq:Poisson} has a unique solution.
% , and we prescribe a Neumann boundary condition.
The derivative of $j$ is given by
\begin{equation}\label{eq:derivative_A2}
    j_u(\xi,u)[v]
    =
    \int_U \left(\beta  |u|^{p-2}u v - q v\right) \dee x,
\end{equation}
where $q \in H^1(U)$ is the adjoint variable and satisfies
\begin{equation*}
    \int_U \left( {\color{black} \kappa(\xi)} \nabla q \cdot \nabla \eta + q \eta + 5 y^4 \eta q \right) \dee x = - \int_U (y(\xi,\cdot)-y_d) \eta \ \dee x,
\end{equation*}
for all $\eta \in H^1(U)$.

Much like in the previous application, we may characterize the steepest descent $v(u)$.
\begin{lem}[Finding the direction of descent]
    Suppose $j_u(\xi,u) \neq 0$.
    One finds that $v(u) = \mu ^{-\frac{1}{(p-1)}}\tilde{v}$, where $\mu = \left(-p^{-1} j_u(\xi,u)[\tilde{v}]\right)^{\frac{p-1}{p}}$ and $\tilde{v}\in \xx$ solves
    \begin{equation}
        p \int_{U} |\tilde{v} |^{p-2} \tilde{v} \eta\ \dee x = - j_u(\xi, u)[\eta] \quad \forall \eta \in \xx.
    \end{equation}
\end{lem}
\begin{proof}
The proof of this follows along the same lines as the proof of Lemma \ref{lem:direct}, thus we omit it.
\end{proof}

\subsubsection{Verification of conditions}
%\todo{Philip: Perhaps rewrite this so that it is not copied and pasted before, maybe make reference to it?}
As mentioned previously, we have made assumptions on $j$, $\mathcal{J}$, and $\xx$ which require verification.
% For the results we have shown, we have made assumptions on $j$, $\mathcal{J}$, and $\xx$.
% These are namely that
For convenience we repeat them here
(i) $\xx$ is the dual of a normed space,
(ii) $j$ is appropriately smooth,
%(iii) $j$ is weak-$*$-lower-semi-continuous in the second variable
(iii) $v \mapsto j_u(\cdot,u)[v]$ is weak-$*$-lower-semi-continuous, and that
%(v) $\mathcal{J}$ is weak-$*$-coercive.
For simplicity, let us assume that $p\geq 2$.
\begin{itemize}
\item[(i)] Since $\xx = L^p(U)$ is reflexive then it is certainly the dual of a space.

\item[(ii)] $j$ is differentiable in the second variable, the form of the derivative is given in \eqref{eq:derivative_A2}
% To see that $j$ is differentiable in the second variable, one may calculate
% \begin{equation}
%     j_u(\cdot,u)[v]
%     =
%     \int_U \left( \alpha |u|^{p-2} u v + (y-y_d) y'(u)[v] \right) \dee x,
% \end{equation}
% where $y'(u)[v] \in H^1(U)$ is the representative of the derivative of $u \mapsto y(u)$ in direction $v \in L^{p}(U)$.
The fact that first derivative is locally Lipschitz follows by seeing that the maps $u \mapsto y(u)$ and $y \mapsto \frac{1}{2}\int_U (y-y_d)^2 \dee x $ are smooth, hence locally Lipschitz in each derivative, and using the same argument as in the preceding section to handle the $p$ power terms.
%\ph{We could do something more, but that might be a lot of a pain!}
%\item[(iii)] For the term $\frac{1}{2} \int_U(y-y_d)^2 \ \dee x$, we have that this is a smooth map composed with $u \mapsto y(\xi,u)$, a compact map which yields continuity.
%We note that we know this is a compact map by standard theory of semi-linear PDEs with monotonoic non-linearity.
%For the remaining term, $\frac{\beta}{p} \int_U |u|^p \dee x$ this is the $p$ power of the norm hence is weak-$*$-lower-semi continuity holds.

\item[(iii)] By reflexivity, this holds.

%\item[(v)] One may show the stronger coercivity assumption by seeing that
%\begin{equation*}
%\end{equation*}
%After integrating in probability,
%\[
%    \mathcal{J}(u) \geq \frac{\beta}{p} \int_U | u|^p\dee x,
%\]
%which yields the desired coercivity.
\end{itemize}

\section{Numerical examples}
\label{sec:num}
This section is devoted to providing numerical simulations to each of the applications presented in Section \ref{sec:app}.
We will test the SSD on each application and for the final application - which fits into a Hilbertian framework, we directly compare the method with the SGD.
We will attempt to numerically verify the condition that $\mathcal{K}_n \leq 0$ as well as measure the decay of $\mathcal{K}_n$, comparing it to the result shown in Theorem \ref{thm:main}.
Throughout this section, we will consider the finite domain $U = [0,1]^2$.
For experiments involving partial differential equations, we will use linear finite elements.
These finite elements are implemented using \texttt{DUNE} \cite{DUNE}, in particular the \texttt{Python} bindings \cite{DUNEPYTHON}.
The sequence of step sizes will be given by $t_n := n^{-1}$ for $n \geq 1$.

\subsubsection{Random coefficient}
\textcolor{black}{Before conducting our experiments we now describe how we simulate our random coefficient for our numerical experiments (excluding Application 1 in Section \ref{subs:app1})}. We take the particular example of defining it on a two-dimensional domain for simplicity. 
We let $-\triangle$ denote the Laplacian, which on the computational domain $U$ is subject to homogeneous Neumann
boundary conditions. From this we can define our covariance operator of our random coefficient $\xi$,
$$C_0=\left(-\triangle+\tau^2\right)^{-\alpha},$$
where $\tau \in \R^{+}$  denotes the inverse lengthscale 
of the random field and $\alpha \in \R^+$ determines the
regularity; specifically draws from the random field are H\"older with
exponent up to $\alpha-1$ (since spatial dimension $d=2$). 
From this we note that the eigenvalue problem
$$C_0\varphi_k=\lambda_k \varphi_k,$$
has solutions, %for $\mathbb{Z}=\{0,1,2,\cdots\}$,
$$\varphi_k(x)=\sqrt{2}\cos(k\pi  x), \quad \lambda_k=\left(|k|^2\pi^2+\tau^2\right)^{-\alpha}, \quad k \in \NN^2.$$
Here $X=L^2(U,\R)$ and the $\varphi_k$ are orthonormal in
$X$ with respect to the standard inner-product.
Draws from a Gaussian measure $\mu_{C_0}$ (Given by Definition \ref{def:gauss}), with covariance $C_0$ are given by the Karhunen-Lo\`{e}ve expansion (KLE)
\begin{equation}
\label{eq:KL}
\Theta = \sum_{k \in \NN^2} \sqrt{\lambda_k}\xi_k \varphi_k(x), \quad \xi_k \sim N(0,1),\quad\rm{i.i.d.}\, .
\end{equation}
{This random function will be almost surely in $X$ and in $C(U,\R)$
provided that $\alpha>d/2$, therefore we 
impose this condition. As we are working in a 2$d$ setting we require $\alpha>1$.
Random draws from the KLE \eqref{eq:KL} are provided in Figure \ref{fig:KL} for varying lengthscale and  regularity.

It is typical in the computational setting to cut off the sum in for the KLE, known as a truncation.
\begin{figure}[h!]
\centering
\includegraphics[width=0.9\textwidth]{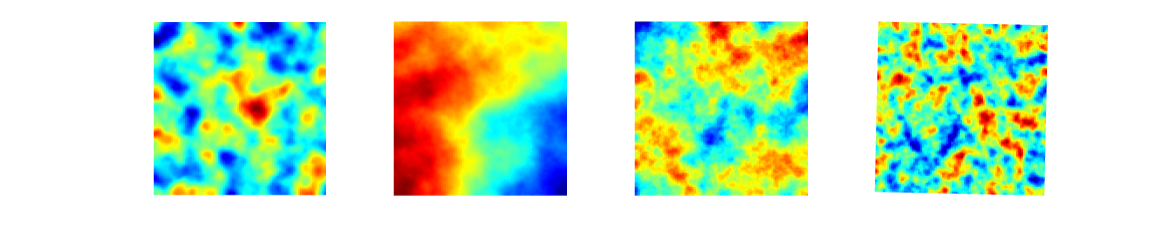}
\caption{Random draws from a Karhunen-Lo\`{e}ve expansion on the square domain $U = [0,1]^2$. Left: $\alpha=2.6, \tau = 25$. Left middle: $\alpha=1.9, \tau = 10$. Right middle: $\alpha=1.7, \tau = 20$. Right: $\alpha=2.3, \tau = 40$}.
\label{fig:KL}
\end{figure}

\subsection{Application 1}\label{sec:exp:basic}
For the first experiment, we consider the minimization of the one-dimensional example which is described in Section \ref{subs:app1}.
We start at an initial guess of $u_0 = -3$ and apply the SSD method.
The results of this experiment - the energy and derivative quantities are provided in Figures \ref{fig:app1:energy} and \ref{fig:app1:derivatives}.
Let us note that this problem does satisfy all of the assumptions for our theorem and does appear to converge, however it does not converge to the minimizer.
One might attribute this to the fact that in $(-1,1)$ it holds that $\mathbb{E}[\mathcal{J}'(u)[v(u)]] = 0$ as in Lemma \ref{lem:OneDSteepestDescent}.
\begin{figure}[h!]
    \centering
    \begin{tikzpicture}
        \begin{axis}
            \addplot table [x = Iteration, y = Energy, col sep = comma]{data_two_point.csv};
            \addlegendentry{$J(u_n)$};
            \addplot table [x = Iteration, y = sampleEnergy, col sep = comma, only marks]{data_two_point.csv};
            \addlegendentry{$j(\xi_n,u_n)$}
        \end{axis}
    \end{tikzpicture}
    \caption{The sample energy and average energy for the 1D toy problem experiment in Section \ref{sec:exp:basic}.}
    \label{fig:app1:energy}
\end{figure}
\begin{figure}[h!]
    \centering
    \begin{tikzpicture}
        \pgfplotsset{unbounded coords=discard}
        \begin{axis}
            \addplot table [x = Iteration, y expr = -\thisrow{Kn}, col sep = comma]{data_two_point.csv};
            \addlegendentry{$-K_n$};
            \addplot table [x = Iteration, y expr = {-\thisrow{KnMax}}, col sep = comma]{data_two_point_processed.csv};
            \addlegendentry{$\min_{i\leq n} (-K_i)$};
            \addplot table [x = Iteration, y expr = -\thisrow{sampleDescent}, col sep = comma, only marks]{data_two_point.csv};
            \addlegendentry{$-J_u(\xi_n,u_n)[v_n]$};
        \end{axis}
    \end{tikzpicture}
    \caption{Derivative information for the 1D toy problem experiment in Section \ref{sec:exp:basic}.}
    \label{fig:app1:derivatives}
\end{figure}

\subsection{Application 2}\label{sec:exp:plaplace}
For this experiment we consider the $p$-Laplace problem described in Section \ref{sec:app:p-laplace}, where we choose $p=4$ and take the random data $g$ to be given by $g = \Theta + \cos(\pi x_1)^2 \cos(\pi x_2)^2$, where $\Theta$ is the KLE with prescribed hyperparameters of the lengthscale $\tau=1$ and the smoothness $\alpha=3$.
For computational simplicity, when we calculate $\Theta$, we cut off so that $k_1 ,k_2 \leq 10$.
We also take $t_k = \frac{1}{10 (k+1)}$.
Our aim from this numerical experiments is to verify the convergence result from our main theorem.

The results of this experiment are provided in Figure \ref{fig:pLaplace:Derivatives}.
%\todo{I think I didn't put out enough data into the csv -- sorry, I should have also included the $-j'(u_h,\xi_h)$}
It is seen that the value $K_n$ appears to abide by the result of the theorem.
%\todo{Philip: can we chat about this.  I would like to give some idea of how good/bad the approximation is.  I can rerun the code calculating the variance of $K_n$ for this setting if that helps.  This was the slow code, but I can make it better by only doing 50 iterations, like the other code, as well as a more coarse grid}

% Figure \ref{fig:pLaplace:Derivatives} provides the negative of the directional derivative along the iterations $n \in \NN$, along with the theoretical rate shown in \eqref{eq:rate}. As we can observe the numerical simulations appear to abide by the rate attained in Theorem \ref{thm:main}.
\begin{figure}[htbp]
    \centering
    \begin{tikzpicture}
        \begin{axis}[ymode = log]
            \addplot table [x = Iteration, y expr = -\thisrow{Kn}, col sep = comma]{data_no_control.csv};\addlegendentry{$-K_n$}
            \addplot table [x = Iteration, y expr = -\thisrow{KnMax}, col sep = comma]{data_no_control_processed.csv};\addlegendentry{$\min_{i \leq n}(-K_i)$}
            \addplot table [x = Iteration, y expr = -\thisrow{sampleDescent}, col sep = comma, only marks]{data_no_control.csv};\addlegendentry{$-j_u(u_n, \xi_n)[v_n]$}
            \addplot table [x = Iteration, y expr = 40/\thisrow{rate}, col sep = comma, dotted, mark = none]{data_no_control_processed.csv};\addlegendentry{$4\left(\sum_{k=1}^n \frac{1}{t_k} \right)^{-1}$}
        \end{axis}
    \end{tikzpicture}
    \caption{Derivative information for the $p$-Laplace-type experiment in Section \ref{sec:exp:plaplace}.}
    \label{fig:pLaplace:Derivatives}
\end{figure}

\subsection{Application 3}\label{sec:exp:4Control}
The next numerical experiment will be the PDE constrained optimization problem described in Section \ref{sec:app:optControl}.
Here we will again choose $p = 4$.
For our random coefficient we choose $D = I_h (1 + \exp(\Theta))$, where $I_h$ is the standard Lagrange interpolation and again $\Theta$ is the KLE with prescribed hyperparameters for the lengthscale $\tau=1$ and smoothness $\alpha=2$. As before, for computational purposes, we again cut off so that $k_1 ,k_2 \leq 10$.
We choose $$y_d(x) = 1 + 256 ((x_1-1)x_1(x_2-1) x_2 )^2,$$ where $\beta = 10^{-2}$, and $F(\xi) = 1 + 5 \Theta$, where this $\Theta$ is the same KLE as above, but drawn independently.
As before we are interested in observing the convergence.
%$n$ the effect of the minimiser as we increase the number of iterations $n$, in terms of verifying the theoretical rate.
The results of this experiment are provided in Figure \ref{fig:4Control:derivative}.
Similarly to the $p$-Laplace example in Section \ref{sec:exp:plaplace}, the value of $K_n$ appears to obey the result of Theorem \ref{thm:main}, but not necessarily the hypothesis.
This is fairly consistent, as the rate matches the directional derivative and minimum for almost all iterations.
\begin{figure}[htbp]
    \centering
    \begin{tikzpicture}
        \begin{axis}[ymode = log]
            \addplot table [x = Iteration, y expr = -\thisrow{Kn}, col sep = comma]{data_pde_control_p4_.csv};\addlegendentry{$-K_n$}
            \addplot table [x = Iteration, y expr = -\thisrow{KnMax}, col sep = comma]{data_pde_control_p4__processed.csv};\addlegendentry{$\min_{i\leq n}(-K_i)$}
            \addplot table [x = Iteration, y expr = -\thisrow{sampleDescent}, col sep = comma, only marks]{data_pde_control_p4_.csv};\addlegendentry{$-j_u(u_n,\xi_n)[v_n]$}
            \addplot table [x = Iteration, y expr = .025/\thisrow{rate}, col sep = comma, dotted, mark = none]{data_pde_control_p4__processed.csv};\addlegendentry{$\frac{1}{40}\left( \sum_{k=1}^n \frac{1}{t_k} \right)^{-1}$}
        \end{axis}
    \end{tikzpicture}
    \caption{Derivative information for the optimal control experiment in Section \ref{sec:exp:4Control}.}
    \label{fig:4Control:derivative}
\end{figure}
\subsection{Application 3\texorpdfstring{$^\prime$}{dash}: A comparison of SGD and SSD}\label{sec:exp:2Control}
Our final numerical experiment is to provide a comparison between both the SSD and SGD, which will hopefully demonstrate the improvements for the former method. 
For this experiment we consider the exact setting of the previous application but instead we set $p=2$, to ensure that we are in a Hilbertian setting. As we are in the Hilbertian setting, one may apply the classical stochastic gradient descent method.
It is easy to verify that in this setting, if one were to consider stochastic steepest descent with step length $S_n := \|J'(u_n)\|_{\xx^*} t_n$, one recovers SGD with step length $t_n$.
Heuristically, one might then expect that gradient descent methods, such as SGD, will work well far away from the minimum, when the derivative has a large magnitude, whereas the SSD should work faster when the derivative has a small magnitude.
For our comparison, we consider the energy, since this is a random problem, it is more realistic to take a Monte-Carlo sample, this is given as
$$
    \frac{\beta}{p}\int_U |u_n|^p \dee x + \frac{1}{2N}\sum^{N}_{i=1}\int_U  (y(\xi_{n_i},u_n) -y_d)^2 \dee x,
$$
where we have $p=2$, and choose $N = 20$ samples. We emphasize the choice of samples is not crucial, where we have repeated the experiment for $N=200$ samples with the same phenomenon occurring. 
% where again we state that the first term is a penalty term.
\begin{figure}
    \centering
    \begin{tikzpicture}
        \begin{axis}[]
            \addplot table [x = Iteration, y = Energy, col sep = comma]{data_pde_control_p2_.csv};\addlegendentry{{\rm SSD}}
            \addplot table [x = Iteration, y = Energy, col sep = comma]{data_unnormalised_pde_control_p2_.csv};\addlegendentry{{\rm SGD}}
        \end{axis}
    \end{tikzpicture}
    \caption{Comparison of the energies for the optimal control experiments in Section \ref{sec:exp:2Control}.}
    \label{fig:2Control:Energy}
\end{figure}
\begin{figure}
    \centering
    \begin{tikzpicture}
        \begin{axis}[ymode = log]
            \addplot table [x = Iteration, y expr = -\thisrow{sampleDescent}, col sep = comma]{data_pde_control_p2_.csv};\addlegendentry{$-j_u(u_n,\xi_n)[v_n]$ (SSD)}
            \addplot table [x = Iteration, y expr = -\thisrow{sampleDescent}, col sep = comma]{data_unnormalised_pde_control_p2_.csv};\addlegendentry{$-j_u(u_n,\xi_n)[v_n]$ (SGD)}
            % K_n
            % -min K_n
            % theoretical rate
        \end{axis}
    \end{tikzpicture}
    \caption{Comparison of the derivatives for the optimal control experiments in Section \ref{sec:exp:2Control}.}
    \label{fig:2Control:derivative}
\end{figure}
Our first comparison is provided in Figure  \ref{fig:2Control:Energy} 
where we compare both methods in terms of the energy for optimal control experiments. As we can observe, immediately after the first iteration, the energy decreases much more quickly and remains the same for the SSD. This is unlike the SGD method whose energy is slower to decrease and remains slightly higher with a difference of around 0.04. Finally we decide to compare both methods in terms of optimal control problem, related to the derivatives. Our results are given in Figure \ref{fig:2Control:derivative}. The figure provides two important insights, firstly the theoretical rate obtained matches the decay of SSD method. Secondly we see after 100 iterations the difference in the derivative is of the order $\mathcal{O}(10^{-2})$, demonstrating an improvement with the SSD.

From Figures \ref{fig:2Control:Energy} and \ref{fig:2Control:derivative},
we comment that there is effectively no additional cost to calculate the steepest descent, as this is only a re-scaling of the gradient.
This can of course be different if one were to have a norm whose derivative does not satisfy an appropriate homogeneity property.

\begin{remark}
We conclude the numerical section with a brief discussion on the choice of $t_n$. As stated previously, the choice we have used is similar to that of a Robbins-Munro condition for stochastic gradient methods. in the SSD setting, the choice is slightly more sensitive, however a to understand this a complexity analysis is required. We leave this for future work.
\end{remark}

% , e.g., $\|u\|_\xx = \|u\|_{L^2} + \|u\|_{L^1}$.
% , which may prove useful to analyze for future work.
}
%\section{What is to be done}
%\textcolor{blue}{
%\begin{itemize}
%\item Add the toy example as a potential easy application to begin with.
%\item Answer the theoretical questions/queries they have. 
%\item Add an in-depth description of the numerical simulations. DONE!!
%\item Get rid of mathbb and use \ instead, i,e, $\R$. Do the same for N %and Z, i.e. use \ .
%%\item Do the graph for $\mathcal{J}$.
%\end{itemize}}

\section{Conclusion}
\label{sec:conc}

The purpose of this paper was to provide a first understanding of the stochastic steepest descent method in a Banach space setting. Commonly gradient methods are exploited which naturally induce a Hilbertian setting. In this work we provided a first simple understanding of a convergence analysis in a general Banach space setting, where one does not require the assumption of reflexivity. Our main result also include a rate of convergence which is similar to that which is achieved for the stochastic gradient method. Numerical simulations were conducted comparing the stochastic steepest descent method to the gradient method, which was tested on two problems: a random 2D elliptic Poisson problem and a $p$-Laplacian problem. Our results demonstrate the improvements of our methodology, related to both relative errors and the energy.

There are numerous natural extensions one can consider from this this work. We provide a short summary of some of these below.
\begin{itemize}
\item One direction would be shape-optimization which provided the initial motivation behind this work.
In particular one could aim to understand recent Banach space algorithms which include the $W^{1,\infty}$-approach, which has been discussed in the various works \cite{DJH22,PJH23}. New convergence analysis is required here.
\item Our numerical experiments were conducted on PDE-constrained optimization problems, which induce randomness. We assumed the randomness to be independent, however in various scenarios one instead has correlated data. Extending this to that setting is of particular interest.
\item Related to the first point, one may consider the use of such a method in parameter estimation problems, or inverse problems. Such problems assume a  unique minimizer exists in a Hilbert space, or use gradient methodologies which relate to that. Such a work would go beyond what has been done in \cite{JK23}, and consider fully the SSD. Other potential connections could lie within bilevel optimization \cite{CST22}, and in the context of non-Gaussian reconstruction \cite{CLR19,SCR22}. 
\item One could consider a more rigorous convergence analysis related to Theorem \ref{thm:main}. The hope is that one could aim to remove the liminf condition. This could be achieved using ideas from gradient domination, which could also verify the convergence of the loss functions. Such results exist for the SGD case in \cite{WKA24}.
\item A final direction for work is related to the above point, but in the context of the step-size. It would be of interest to find, an optimal choice of $t_n$, using the traditional analysis done for SGD, and related methods. This exhaustive work would result in a new paper, and thus we leave this for future work.
\end{itemize}

\section*{Acknowledgements}
NKC is supported by an EPSRC-UKRI AI for Net Zero Grant: ``{Enabling CO2 Capture And Storage Projects Using AI}", (grant EP/Y006143/1). NKC is also supported by a City University of Hong Kong Start-up Grant, project number 7200809
PJH was funded by EPSRC (grant EP/W005840/1) during the preparation of this work.

\end{document}